\documentclass[11pt]{article}
\usepackage{amsmath}
\usepackage{mathrsfs}
\usepackage{amsthm}
\usepackage{amsfonts}
\usepackage{amssymb}
\usepackage{bm}
\usepackage{mathtools}
\usepackage{color}
\usepackage{tikz}

\usepackage{enumerate}
\usepackage{enumitem}
\setlist[enumerate,1]{label=(\arabic*),font=\textup,
leftmargin=7mm,labelsep=1.5mm,topsep=0mm,itemsep=-0.8mm}
\setlist[enumerate,2]{label=(\alph*).,font=\textup,
leftmargin=7mm,labelsep=1.5mm,topsep=-0.8mm,itemsep=-0.8mm}

\usepackage[pdfstartview=FitH,bookmarksopen=false,
colorlinks=true,linkcolor=blue,citecolor=blue]{hyperref}

\newtheorem{theorem}{Theorem}[section]

\newtheorem{corollary}{Corollary}[section]
\newtheorem{lemma}{Lemma}[section]

\title{\bf Extremal problems for disjoint graphs \thanks {Research was partially supported by the National
Nature Science Foundation of China (grant numbers 11871329, 11971298, 12201161), Hainan Provincial Natural Science Foundation of China (No. 122QN218),  Doctoral Scientific Research Foundation of Hainan University (No. KYQD(ZD)-21153) and Doctoral Scientific Research Foundation of Henan Normal University (No. QD2022053).} }
\author {Zhenyu Ni$^{1}$, \,Jing Wang$^{2}$,  \,  Liying Kang$^{3}$\thanks{\em Corresponding author. Email address: lykang@shu.edu.cn (L. Kang), 1051466287@qq.com(Z. Ni), wj517062214@163.com(J. Wang).}\\
{\small $^{1}$Department of Mathematics, Hainan University,
Haikou 570228, P.R. China}\\
{\small $^{2}$College of Mathematics and Information Science, Henan Normal University,}
\\{\small Xinxiang 453007, P.R. China}\\
{\small $^{3}$Department of Mathematics, Shanghai University,
Shanghai 200444, P.R. China}}

\date{}
\begin{document}

\maketitle

\begin{abstract}

For a simple graph $F$, let
$\mathrm{EX}(n, F)$ and $\mathrm{EX_{sp}}(n,F)$ be  the  set of graphs with the maximum number of edges and the set of graphs with the maximum spectral radius in an $n$-vertex graph without any copy of the graph $F$, respectively.
Let $F$ be a graph with $\mathrm{ex}(n,F)=e(T_{n,r})+O(1)$.
 In this paper, we show that $\mathrm{EX_{sp}}(n,kF)\subseteq \mathrm{EX}(n,kF)$ for sufficiently large $n$. This generalizes a result of Wang, Kang and Xue [J. Comb. Theory, Ser. B, 159(2023) 20-41]. We also determine the extremal graphs of $kF$ in term of the extremal graphs of $F$.

\bigskip \noindent{\bf Keywords:} Extremal problem; Tur\'{a}n number; Spectral radius; Spectral extremal problem

\medskip

\noindent{\bf AMS (2000) subject classification:}  05C50; 05C35
\end{abstract}

\section{Introduction}
Given a family of graphs $\mathcal F$, a graph $G$ is called $\mathcal F$-{\sl free} if it does not contain any copy of $F\in \mathcal F$. If $\mathcal F=\{F\}$, we simply use $F$-free. 
The {\sl Tur\'{a}n number} of $\mathcal F$, denoted by $\mathrm{ex}(n,\mathcal F)$, is the maximum number of edges among all $n$-vertex $\mathcal F$-free graphs. Denote by $\mathrm{EX}(n,\mathcal F)$ the set of all $n$-vertex $\mathcal F$-free graphs with $\mathrm{ex}(n,\mathcal F)$ edges. The problem of determining 
$\mathrm{ex}(n,\mathcal F)$ and $\mathrm{EX}(n,\mathcal F)$ is usually called the Tur\'{a}n type extremal problem. 

Let $kF$ be the disjoint union of $k$ copies of $F$. The {\sl join} of $G$ and $H$, denoted by $G\vee H$, is obtained from a disjoint union of $G$ and $H$ by adding edges joining every vertex of $G$ to vertex of $H$. Let $K_{r}$ denote the complete graph of order $r$.
The Tur\'an graph $T_{n,r}$ is the complete $r$-partite graph on $n$ vertices with its part sizes are as equal as possible.

In 1941, Tur\'{a}n \cite{Turan1941} proved that if $G$ is $K_{r+1}$-free of order $n$, then $e(G)\leq e(T_{n,r})$, with equality if and only if $G=T_{n,r}$. Simonovits \cite{Simonovits1968} and Moon \cite{JWMoon1968} showed that if $n$ is sufficiently large, $K_{k-1}\vee T_{n-k+1}$ is the unique extremal graph for $kK_{r+1}$. In 2011, Gorgol \cite{Gorgol2011} considered the Tur\'{a}n number of disjoint copies of connected graphs, and prove that $\mathrm{ex}(n,kF)=\mathrm{ex}(n,F)+O(n)$. Yuan and Zhang \cite{YuanZhang2017} determined the Tur\'{a}n number and the extremal graphs of disjoint copies of $P_3$. For recent results about Tur\'{a}n numbers of disjoint copies of graphs see \cite{BushawKettle2011,GerbnerMethukuVizer2019,LiuLidickyPalmer2013}.



The {\sl spectral radius} of a graph $G$ is the spectral radius of its adjacent matrix $A(G)$, denoted by $\rho(G)$.
In 2010, Nikiforov \cite{Nikiforov2010extremalspectra} proposed the spectral Tur\'{a}n type extremal problem which determines the maximum spectral radius among all $n$-vertex $\mathcal F$-free graphs. Let $\mathrm{EX}_{sp}(n,\mathcal F)$ denote the set of graphs which have the maximum spectral radius among all $n$-vertex $\mathcal F$-free graphs. Several results are known in this area (see \cite{CDT21,Yongtao21,DKLNTW,LZZ2022,zhai2022}).

In 2007, Nikiforov \cite{Nikiforov07} determined the spectral version of Tur\'{a}n's extremal result, and proved that if $G$ is a $K_{r+1}$-free graph on $n$ vertices, then $\rho (G)\le \rho (T_{n,r})$, with equality if and only if $G=T_{n,r}$. Cioab\u{a}, Feng, Tait and Zhang \cite{CioabaFengTaitZhang} proved that the graph attaining the maximum spectral radius among all $n$-vertex $F_k$-free graphs is a member of $\mathrm{EX}(n,F_k)$, where $F_k$ is the graph consisting of $k$ triangles which intersect in exactly one common vertex. In 2022, Cioab\u{a}, Desai and Tait \cite{CDT21} raised the following conjecture: Let $F$ be any graph such that the graphs in $\mathrm{EX}(n,F)$ are Tur\'{a}n graphs plus $O(1)$ edges. Then for sufficiently large $n$, $\mathrm{EX_{sp}}(n, F)\subseteq \mathrm{EX}(n, F)$.
In 2023, Wang, Kang and Xue \cite{Kang2023} proved the Cioab\u{a}-Desai-Tait conjecture, and gave a stronger result: Let $F$ be any graph such that $\mathrm{ex}(n,F)=e(T_{n,r})+O(1)$. Then $\mathrm{EX_{sp}}(n, F)\subseteq \mathrm{EX}(n, F)$ for sufficiently large $n$. Recently, Ni, Wang and Kang \cite{NiWangKang} proved that
 $\mathrm{EX_{sp}}(n, kK_{r+1})\subseteq \mathrm{EX}(n,kK_{r+1})$ for sufficiently large $n$.

For a graph $F$ with $\mathrm{ex}(n,F)=e(T_{n,r})+O(1)$,
we prove the following theorems for $kF$.

\begin{theorem} \label{edge extrema}
Let $r\geq 2$, $k\geq 1$ be integers, and $F$ be any graph such that $\mathrm{ex}(n,F)=e(T_{n,r})+O(1)$. For sufficiently large $n$, if $G$ has the maximal edges over all $n$-vertex $kF$-free graphs, then
$G = K_{k-1}\vee G(n-k+1,F)$, where $G(n-k+1,F)\in \mathrm{EX}(n-k+1,F)$.
\end{theorem}

\begin{theorem}\label{spectral extrema}
Let $r\geq 2$, $k\geq 1$ be integers, and $F$ be any graph such that $\mathrm{ex}(n,F)=e(T_{n,r})+O(1)$. For sufficiently large $n$, if $H$ has the maximal spectral radius over all $n$-vertex $kF$-free graphs, then
$$H \in \mathrm{EX}(n, kF).$$
\end{theorem}

By Theorems \ref{edge extrema} and \ref{spectral extrema}, we can easily get the following results.

\begin{corollary}[\cite{Kang2023}]
Let $r\geq 2$ be an integer, and $F$ be a graph with $\mathrm{ex}(n, F)=e(T_{n,r})+O(1)$.
 For sufficiently large $n$, if $G$ has the maximal spectral radius over all $n$-vertex $F$-free graphs, then
$$G \in \mathrm{EX}(n, F).$$
\end{corollary}

\begin{corollary}[\cite{NiWangKang}]
For $k\geq 2$, $r\geq 2$, and sufficiently large $n$. Suppose that $G$ has the maximum spectral radius among all $kK_{r+1}$-free graphs on $n$ vertices, then $G$ is isomorphic to  $K_{k-1}\vee T_{n-k+1,r}$.
\end{corollary}

\section{Preliminaries}\label{sec2}

Given a simple graph $G=(V(G),E(G))$, we denote by $\delta(G)$ and $\Delta(G)$ the minimum and maximum degrees of $G$. 
For $V_1,V_2 \subseteq V(G)$, we denote by $e_G(V_1,V_2)$ the number of edges of $G$ with one end in $V_1$ and the other in $V_2$.
For any $v\in V(G)$ and $S\subseteq V(G)$, let $N_S(v)=\{u\in S:~ uv\in E(G)\}$ and $d_{S}(v)=|N_{S}(v)|$. We denote by $G\setminus S$ the graph obtained from $G$ by deleting all vertices of $S$ and all the edges incident with $S$. Let $G[S]$ denote the graph induced by $S$ whose vertex set is $S$ and whose edge set consists of all edges of $G$ which have both ends in $S$. A spanning subgraph of $G$ is a subgraph whose vertex set is $V(G)$.

A set $M$ of disjoint edges of $G$ is called a {\sl matching} in $G$.  The {\sl matching number}, denoted by $\nu(G)$, is the maximum cardinality of a matching in $G$. Let $G$ be a simple graph with matching number $\nu(G)$ and maximum degree $\Delta(G)$. For two given integers $\nu$ and $\Delta$, define $f(\nu, \Delta)
=\max\{e(G): \nu(G)\leq \nu, \Delta(G)\leq \Delta \}$. In 1976, Chv\'atal and Hanson \cite{Chvatal76} obtained the following result.

\begin{theorem}[\cite{Chvatal76}]\label{ffunction}
For integers $\nu \geq 1$ and $\Delta \geq 1$, we have
$$f(\nu, \Delta)= \Delta \nu +\left\lfloor\frac{\Delta}{2}\right\rfloor
 \left \lfloor \frac{\nu}{\lceil{\Delta}/{2}\rceil }\right \rfloor
 \leq \Delta \nu+\nu.$$
\end{theorem}

A graph is $k$-colourable if it has a $k$-colouring such that no two adjacent vertices are assigned the same colour. The minimum $k$ for which a graph $G$ is $k$-colourable is called its chromatic number and denoted $\chi(G)$.
The stability theorem proved by Erd\H{o}s \cite{Erdos1968,Erdos1966} and Simonovits \cite{Simonovits1968} is an important tool in Tur\'{a}n type extremal problem.
\begin{theorem}[\cite{Erdos1968,Erdos1966,Simonovits1968}]\label{edge stability}
For every $\varepsilon >0$, and any graph $F$ with $\chi(F)=r+1\geq 3$, there exists $\delta>0$ such that if a graph $G$ of order $n$ satisfies $e(G)> (1-\frac{1}{r}-\delta)\frac{n^2}{2}$, then either $G$ contains $F$, or $G$ differs from $T_{n,r}$ in at most $\varepsilon n^2$ edges.
\end{theorem}




 Let $A(G)$ be the adjacency matrix of a graph $G$. By the Perron-Frobenius theorem, the spectral radius $\rho(G)$ is the largest eigenvalue of $A(G)$.
Let $\mathbf{x}=(x_1,\cdots,x_n)^{\mathrm{T}}$ be an eigenvector corresponding to $\rho(G)$. Then for any $ i\in [n]$ and $x_i\neq 0$,
\begin{equation}\label{eigenequation}
\rho(G)x_i=\sum_{ij\in E(G)}x_j.
\end{equation}
By the Rayleigh quotient, we have
\begin{equation}\label{Rayleigh}
\rho(G)=\max_{\mathbf{x}\in \mathbb{R}^{n}_{+}}\frac{\mathbf{x}^{\mathrm{T}}A(G)\mathbf{x}}{\mathbf{x}^{\mathrm{T}}\mathbf{x}}=\max_{\mathbf{x}\in \mathbb{R}^{n}_{+}}\frac{2\sum_{ij\in E(G)}x_ix_j}{\mathbf{x}^{\mathrm{T}}\mathbf{x}}.
\end{equation}

%
For a connected graph $G$, the eigenvectors corresponding to $\rho(G)$ are positive, and $G$ is connected if and only if $A(G)$ is irreducible. By \cite{Zhan13}, we have the following result.

\begin{lemma}[\cite{Zhan13}]\label{subgraph}
Let $G$ be a connected graph. If $G'$ is a proper subgraph of $G$, then $\rho(G')< \rho(G)$.
\end{lemma}

Nikiforov \cite{Niki09JGT} proved the spectral version of the stability theorem, the spectral stability theorem can be written as follows.


\begin{theorem}[\cite{Niki09JGT}] \label{stability}
Let $F$ be a graph with chromatic number $\chi (F)=r+1$. For every $\varepsilon >0$, there exist $\delta >0$ and $n_0$ such that if  $G$ is an $F$-free graph on $n\ge n_0$ vertices  with $\rho (G) \ge (1- \frac{1}{r} -\delta )n$, then $G$ can be obtained from $T_{n,r}$ by adding and deleting at most $\varepsilon n^2$ edges.
\end{theorem}

The following lemma was given in \cite{CioabaFengTaitZhang}.
\begin{lemma}[\cite{CioabaFengTaitZhang}]
\label{intersect}
Let $V_1,\cdots,V_k$ be $k$ finite sets. Then
\[
|V_1 \cap \cdots \cap V_k| \geq \sum_{i=1}^{k}|V_i|-(k-1)\left|\bigcup_{i=1}^{k}V_i\right|
\]
\end{lemma}

%

\section{Proof of Theorem \ref{edge extrema}}\label{edge}
In the rest of this paper, we always assume that $F$ satisfies $\mathrm{ex}(n,F)=e(T_{n,r})+a$, and $G$ has the maximum number of edges among all $n$-vertex $kF$-free graphs. By the Erd\H{o}s-Stone-Simonovits theorem \cite{ES66,ES46}, we have
 \begin{equation*}
  \mathrm{ex}(n,F) = \left( 1- \frac{1}{\chi (F) -1}  \right)
  \frac{n^2}{2} + o(n^2).
  \end{equation*}
Combining with the fact $(1-\frac{1}{r})\frac{n^2}{2}-\frac{r}{8}\leq e(T_{n,r})\leq (1-\frac{1}{r})\frac{n^2}{2}$, we have $\chi(F)=r+1$.

Let $G(n-k+1,F)\in \mathrm{EX}(n-k+1,F)$. Then $e(G(n-k+1,F))=e(T_{n-k+1,r})+a$. The aim of this section is to prove that $G$ is the join of $K_{k-1}$ and $G(n-k+1,F)$ for sufficiently large $n$.

\begin{lemma}\label{construct lower bound of edge extrema}
Let $G$ be an $n$-vertex $kF$-free graph with the maximum number of edges. Then
$$e(G)\geq e(T_{n-k+1,r})+(k-1)n+a-\frac{k}{2}(k-1).$$
\end{lemma}
\noindent{\bfseries Proof.}
We first claim that $K_{k-1}\vee G(n-k+1,F)$ is $kF$-free. Suppose that $K_{k-1}\vee G(n-k+1,F)$ has at most $t$ disjoint copies of $F$. Since $G(n-k+1,F)$ is $F$-free, we have $t\leq k-1$, i.e. $K_{k-1}\vee G(n-k+1,F)$ is $kF$-free.
Then we  have
\begin{align*}
e(G)&\geq e(K_{k-1}\vee G(n-k+1,F))\\
&= \binom{k-1}{2}+(k-1)(n-k+1)+\mathrm{ex}(n-k+1,F)\\
&= (k-1)n-\frac{k}{2}(k-1)+e(T_{n-k+1,r})+a.
\end{align*}
\qed


By Lemmas \ref{edge stability}, \ref{construct lower bound of edge extrema} and similar discussion as in
\cite{Kang2023}, we have the following lemma.

\begin{lemma}\label{partition of edge extrema}
Let $G$ be an $n$-vertex $kF$-free graph with the maximum number of edges.
For sufficiently large $n$, there exists $\varepsilon_1 >0$ such that $V(G)$ has a partition $\{U_1, U_2, \cdots, U_{r}\}$ with $\sum_{1\leq i<j\leq r}e_G(U_i,U_j)$ attaining the maximum, and for any $i\in [r]$
$$\frac{n}{r}-\varepsilon_1 n< |U_i|< \frac{n}{r}+\varepsilon_1 n.$$
Furthermore, for sufficiently small constant $\theta$, let
\[
W:=\bigcup_{i=1}^{r}\big\{v\in U_i: d_{U_i}(v)\geq 2 \theta n\big\}.
\]
Then $|W|\leq \theta n$.
\end{lemma}

In the following, we assume that $\ell$ is an integer that is much larger than $\max\{a,|V(F)|\}$. Suppose that $\varepsilon_2$ is a sufficiently small constant with $\varepsilon_2\ll \theta$. Let
\[
L:=\Big\{v\in V(G): d_G(v)\leq \Big(1-\frac{1}{r}-\varepsilon_2\Big)n\Big\}.
\]
\begin{lemma}\label{WL structure}
If $|L|\ll \varepsilon_2 n$, then \\
$(1)$ For any $u\in W\setminus L$, $G\setminus (W\cup L)$ contains a copy of $T_{kr\ell,r}$ such that $u$ is adjacent to all vertices of the copy of $T_{kr\ell,r}$.\\
$(2)$ For $i\in [r]$ and an arbitrary subset $S\subseteq U_{i}\setminus (W\cup L)$,  $G\setminus (W\cup L)$ contains a subgraph $G'$ such that $G'$ has a spanning subgraph $T_{r|S|,r}$ with $V(G')\cap U_{i}=S$.
\end{lemma}
\noindent{\bfseries Proof.}
$(1)$ For any $u\in W\setminus L$, we have $d(u)>(1-\frac{1}{r}-\varepsilon_2)n$. Without loss of generality, we may assume that $u\in U_1$. Then $d_{U_1}(u)\geq 2\theta n$. By Lemma \ref{partition of edge extrema}, we have
\begin{eqnarray*}
d_{U_1\setminus (W\cup L)}(u)&\geq& d_{U_1}(u)-|W\cup L|\\
&\geq & 2\theta n-\theta n-\varepsilon_2 n\\
&> & k\ell.
\end{eqnarray*}
Let $u_{1,1},\cdots, u_{1,k\ell}$ be the neighbors of $u$ in $U_1\setminus (W\cup L)$. Since $\{U_1, U_2, \ldots, U_{r}\}$ is the  partition of $V(G)$ that  maximizes the number of crossing edges of $G$,  $d_{U_1}(u)\leq \frac{1}{r}d(u)$. For any $2\leq j\leq r$,
\begin{eqnarray}
d_{U_j}(u)&\geq& d(u)-d_{U_1}(u)-(r-2)\Big(\frac{n}{r}+\varepsilon_1 n\Big)\nonumber\\[2mm]
&> &\frac{r-1}{r}\Big(1-\frac{1}{r}-\varepsilon_2\Big)n-(r-2)\Big(\frac{n}{r}+\varepsilon_1 n\Big)\nonumber\\[2mm]
&>& \frac{n}{r^2}-o(n).\label{eqWminusL}
\end{eqnarray}
For any $i\in [r]$ and $v\in U_i \setminus (W\cup L)$, we have $d(v)>(1-\frac{1}{r}-\varepsilon_2)n$ and $d_{U_i}(v)<2\theta n$. Then for any $j\in [r]$ and $j\neq i$, we have
\begin{eqnarray}
d_{U_j}(v)&\geq& d(v)-d_{U_i}(v)-(r-2)\Big(\frac{n}{r}+\varepsilon_1 n\Big)\nonumber\\[2mm]
&> &\Big(1-\frac{1}{r}-\varepsilon_2\Big)n-2\theta n-(r-2)\Big(\frac{n}{r}+\varepsilon_1 n\Big)\nonumber\\[2mm]
&> & \frac{n}{r}-o(n).\label{eqWcupL}
\end{eqnarray}
We consider the common neighbors of $u,u_{1,1},\cdots,u_{1,k\ell}$ in $U_2\setminus (W\cup L)$.
Combining with (\ref{eqWminusL}), (\ref{eqWcupL}) and  Lemma \ref{intersect},  we have
\begin{eqnarray*}
& &|N_{U_2}(u)\cap (\cap_{i\in[k\ell]} N_{U_2}(u_{1,i}))\setminus (W\cup L)|\\[2mm]
&\geq & d_{U_2}(u)+\sum_{i=1}^{k\ell}d_{U_2}(u_{1,i})-k\ell|U_2|-|W|-|L|\\[2mm]
&> &\frac{n}{r^2}-o(n) + k\ell\left(\frac{n}{r}-o(n)\right)- k\ell\left(\frac{n}{r}+\varepsilon_1 n\right)-\theta n-\varepsilon_2 n\\[2mm]
&> & \frac{n}{r^2}-o(n)>k\ell.
\end{eqnarray*}
Let $u_{2,1},\cdots, u_{2,k\ell}$ be the common neighbors of $\{u,u_{1,1},\cdots,u_{1,k\ell}\}$ in $U_2\setminus (W\cup L)$. For an integer $2\leq s\leq r-1$, suppose that $u_{s,1},\cdots,u_{s,k\ell}$ are the common neighbors of $\{u,u_{i,1},\cdots,u_{i,k\ell}: 1\leq i\leq s-1\}$ in $ U_{s}\setminus (W\cup L)$. We next consider the common neighbors of $\{u,u_{i,1},\cdots,u_{i,k\ell}: 1\leq i\leq s\}$ in $U_{s+1}\setminus (W\cup L)$.
By  (\ref{eqWminusL}), (\ref{eqWcupL}) and   Lemma \ref{intersect}, we have
\begin{eqnarray*}
& &|N_{U_{s+1}}(u)\cap (\cap_{i\in [s], j\in [k\ell]} N_{U_{s+1}}(u_{i,j}))\setminus (W\cup L)|\\[2mm]
&\geq & d_{U_{s+1}}(u)+\sum_{i=1}^{s}\sum_{j=1}^{k\ell}d_{U_{s+1}}(u_{i,j})-sk\ell|U_{s+1}|-|W|-|L|\\[2mm]
&> & \frac{n}{r^2}-o(n)+sk\ell\left(\frac{n}{r}-o(n)\right)-sk\ell\left(\frac{n}{r}+\varepsilon_1 n\right)-\theta n-\varepsilon_2 n\\[2mm]
&> & \frac{n}{r^2}-o(n)> k\ell.
\end{eqnarray*}
Let $u_{s+1,1},\cdots, u_{s+1,k\ell}$ be the common neighbors of $\{u,u_{i,1},\cdots,u_{i,k\ell}: 1\leq i\leq s\}$ in $U_{s+1}\setminus (W\cup L)$. Therefore, for every $i\in [r]$, there exist $u_{i,1},\cdots,u_{i,k\ell}\in U_{i}\setminus (W\cup L)$ such that $\{u_{1,1},\cdots,u_{1,k\ell}\}$ $\cup ~\{u_{2,1},\cdots,u_{2,k\ell}\}$ $\cup ~\cdots ~\cup$ $\{u_{r,1},\cdots,u_{r,k\ell}\}$ form a complete $r$-partite subgraph of $G\setminus (W\cup L)$ and $u$ is adjacent to all the above $kr\ell$ vertices.

$(2)$ Without loss of generality, take an arbitrary subset $S\subseteq U_{1}\setminus (W\cup L)$. We consider the common neighbors of $S$ in $U_2\setminus (W\cup L)$. By using the similar discussion as in $(1)$ of Lemma \ref{WL structure}, we can find a copy of $T_{r|S|,r}$ in  $G\setminus (W\cup L)$. The result follows.
\qed

\begin{lemma}\label{WL property}
If $|L|\ll \varepsilon_2 n$, then $|W\setminus L|\leq k-1$, and $d(u)=n-1$ for any $u\in W\setminus L$.

\end{lemma}
\noindent{\bfseries Proof.}
We first prove that $|W\setminus L|\leq k-1$. Suppose to the contrary that $|W\setminus L|\geq k$. By Lemma \ref{WL structure}, for any $u\in W\setminus L$, we can find a copy of $K_1 \vee T_{kr\ell,r}$, say $G_1$, such that $u$ is the vertex with degree $kr\ell$ in $G_1$. Therefore, $e(G_1)=e(T_{kr\ell+1,r})+k\ell>e(T_{kr\ell+1,r})+a=\mathrm{ex}(kr\ell+1,F)$. As $T_{kr\ell+1,r}$ is $F$-free, we can find a copy of $F$ which contains $u$ in $G$. Since $|W\setminus L|\geq k$ and $\ell$ is much larger than $\max\{a,|V(F)|\}$, we can find $k$ disjoint copies of  $F$. This contradicts the fact that $G$ is $kF$-free.

Now we prove that each vertex of $W\setminus L$ has degree $n-1$. Suppose to the contrary that there exists $u\in W\setminus L$ such that $d(u)<n-1$. Assume that $v$ is a vertex of $G$ such that $uv\notin E(G)$. Let $G'$ be the graph with $V(G')=V(G)$ and $E(G')=E(G)\cup \{uv\}$. We claim that $G'$ is $kF$-free. Otherwise, $G'$ contains a copy of $kF$, say $G_2$, then $uv\in E(G_2)$. Therefore, $G_2$ contains a copy of $F$, say $G_3$, such that $uv\in E(G_3)$. Let $G_4=G_2\setminus G_3$. Then $G_4$ is a copy of $(k-1)F$ in $G$ and $u\notin V(G_4)$. Since $u\in W\setminus L$, by Lemma \ref{WL structure}, we can find a copy of $K_1 \vee T_{kr\ell,r}$, say $G_5$,  such that $u$ is the vertex with degree $kr\ell$ in $G_5$. Since $e(G_5)=e(T_{kr\ell+1,r})+k\ell>e(T_{kr\ell+1,r})+a$, there is a copy of $F$, say $G_6$, such that $V(G_6)\cap V(G_4)=\emptyset$. Then $G_4\cup G_6$ is a copy of $kF$ in $G$, and this is a contradiction. Therefore, $G'$ is $kF$-free. By the construction of $G'$, we have $e(G')>e(G)$, which contradicts the assumption that $G$ has the maximum number of edges among all $n$-vertex $kF$-free graphs.
\qed

\begin{lemma}\label{WL edge}
If $|L|\ll \varepsilon_2 n$, then $e(G[U_i\setminus (W\cup L)])<k\ell(k\ell+1)$ for any $i\in [r]$.

\end{lemma}
\noindent{\bfseries Proof.}
We first prove that $\Delta(G[U_i\setminus (W\cup L)])<k\ell$ for any $i\in [r]$. Suppose to the contrary that there exist $i_0\in [r]$ and  $u\in U_{i_0}\setminus (W\cup L)$ such that $d_{U_{i_0}\setminus (W\cup L)}(u)\geq k\ell$. As $u\in U_{i_0}\setminus (W\cup L)$, there exists a vertex $w$ such  that $uw\notin E(G)$. Let $G'$ be the graph with $V(G')=V(G)$ and $E(G')=E(G)\cup \{uw\}$. Then $G$ is a proper subgraph of $G'$. By the maximality of $e(G)$, $G'$ contains a copy of $kF$, say $G_1$. From the construction of $G'$, we see that $u\in V(G_1)$, and $G_1$ contains  a copy of $(k-1)F$, say $G_2$, such that $u\notin V(G_2)$. Obviously, $G_2\subseteq G$. Since $d_{U_{i_0}\setminus (W\cup L)}(u)\geq k\ell$, let $S=N_{U_{i_0}\setminus (W\cup L)}(u)\cup\{u\}$. Then by Lemma \ref{WL structure}, $G$ has a subgraph $G_3$ which contains a spanning subgraph $T_{r(k\ell+1),r}$ such that $e(G_3)\geq e(T_{r(k\ell+1),r})+k\ell$. Since $u\notin V(G_2)$, $G_3$ has a subgraph $G_4$ which contains a spanning subgraph $T_{r\ell,r}$ such that $V(G_4)\cap V(G_2)=\emptyset$ and $e(G_4)> e(T_{r\ell,r})+a$. Therefore, $G_4$ has a copy of $F$ such that this $F$ together with $G_2$ forms a copy of $kF$ in $G$, which contradicts the fact that $G$ is $kF$-free.

Now we claim that for any $i\in [r]$, $\nu(G[U_i\setminus (W\cup L)])< k\ell$. Otherwise, there exists $i_0\in [r]$ such that $G[U_{i_0}\setminus (W\cup L)]$ contains a matching $M$ with $k\ell$ edges. By Lemma \ref{WL structure}, $G$ has a subgraph $G_5$ which contains a spanning subgraph $T_{2kr\ell,r}$ such that with $V(G_1)\cap U_{i_0}=V(M)$ and $e(G_5)\geq e(T_{2kr\ell,r})+k\ell$. If we partition $M$ into $k$ disjoint parts such that each part has $\ell$  matching edges, then we can find $k$ disjoint subgraphs in $G_5$ such that each subgraph has at least $e(T_{2r\ell,r})+\ell$ edges. Therefore, we can find $k$ disjoint copies of $F$ in $G$, a contradiction.

Combining with Lemma \ref{ffunction}, for any $i\in [r]$, we have
\begin{align*}
e(G[U_i\setminus (W\cup L)])&\leq f(\nu(G[U_i\setminus W]),\Delta(G[U_i\setminus W]))\\
&< f(k\ell,k\ell)\\
&\leq k\ell(k\ell+1).
\end{align*} \qed
\begin{lemma}\label{L is empty in edge extrema}
If $L=\emptyset$, then $G= K_{k-1}\vee G(n-k+1,F)$.
\end{lemma}
\noindent{\bfseries Proof.}
We first prove that $|W|=k-1$ if $L=\emptyset$. Otherwise, $|W|<k-1$ by Lemma \ref{WL property}. Let $|W|=s$. By Lemma \ref{WL edge}, for any $i\in [r]$, $e(G[U_i\setminus W])<k\ell(k\ell+1)$. Then $e(\cup_{i=1}^{r}G[U_i\setminus W])<rk\ell(k\ell+1)$.
Take $S\subseteq U_1\setminus W$ with $|S|=k-s-1$. Let $G'$ be the graph with $V(G')=V(G)$ and $E(G')=E(G)\setminus \{uv: uv\in \cup_{i=1}^{r}E(G[U_i\setminus W])\}\cup \{uv: u\in S,v\in (U_1\setminus W)\setminus S\}$. It is obvious that $G'$ is $kF$-free. However, for sufficiently large $n$,
\begin{align*}
e(G')-e(G)
& = |S|(|U_1|-|W|-|S|)- e(\cup_{i=1}^{r}G[U_i\setminus W]) \\[2mm]
& > (k-s-1)\left(\frac{n}{r}-o(n)\right)-rk\ell(k\ell+1)\\[2mm]
&>0.
\end{align*}
This contradicts the fact that $G$ has the maximum number of edges over all $n$-vertex $kF$-free graphs.
Therefore, $|W|=s=k-1$.

By Lemma \ref{WL property}, if $L=\emptyset$, then each vertex in $W$ has degree $n-1$. Since $|W|=k-1$, $G$ is the join of $K_{k-1}$ and a $r$-partite graph of order $n-k+1$. It is sufficient to prove that $G\setminus W\in \mathrm{EX}(n-k+1,F)$. Now we prove that $G\setminus W$ is $F$-free. Otherwise, $G\setminus W$ contains a copy of $F$, say $G_1$. Since $|W|=k-1$ and each vertex in $W$ has degree $n-1$, by Lemma \ref{WL structure} and similar discussion as in Lemma \ref{WL property}, we can find $k-1$ disjoint copies of $F$ which are disjoint with $G_1$ . This contradicts the fact that $G$ is $kF$-free. Therefore, by the structure of $G$ and the maximality of $e(G)$, we have $G\setminus W\in \mathrm{EX}(n-k+1,F)$. The proof of Lemma \ref{L is empty in edge extrema} is completed.
\qed

\medskip
\noindent{\bfseries Proof of Theorem \ref{edge extrema}.} By Lemma \ref{L is empty in edge extrema}, if $L=\emptyset$, there is nothing to do. Now we consider the case that $|L|>0$. Let $G_0=G$ and $L_0=L$. Take $u_1\in L_0$. Let $G_1=G_0\setminus \{u_1\}$ and $L_1=\{v\in V(G_1): d_{G_1}(v)\leq (1-\frac{1}{r}-\varepsilon_2)(n-1)\}$. 
If $|L_1|>0$, take $u_2\in L_1$. Let $G_2=G_1\setminus \{u_2\}$ and $L_2=\{v\in V(G_2): d_{G_2}(v)\leq (1-\frac{1}{r}-\varepsilon_2)(n-2)\}$. 
Continue this process as long as after $t$ steps we get a subgraph $G_t$ and $|L_t|=0$.
For any $i\in [t]$, $G_{i}$ is obtained from $G_{i-1}$ by deleting a vertex with degree small than $(1-\frac{1}{r}-\varepsilon_2)(n-i+1)$ in $G_{i-1}$. By Lemma \ref{construct lower bound of edge extrema}, we have $$e(G_0)\geq e(T_{n-k+1,r})+(k-1)n+a-\frac{k}{2}(k-1)>e(T_{n,r}).$$
Since $\delta(T_{n,r})>\big(1-\frac{1}{r}-\varepsilon_2\big)n,$ we have
\begin{align*}
e(G_1)&\geq e(G_0)-\Big(1-\frac{1}{r}-\varepsilon_2\Big)n\\[2mm]
&> e(T_{n,r})-\Big(1-\frac{1}{r}-\varepsilon_2\Big)n\\[2mm]
&\geq  e(T_{n-1,r})+\delta(T_{n,r})-\Big(1-\frac{1}{r}-\varepsilon_2\Big)n\\[2mm]
&\geq e(T_{n-1,r})+1.
\end{align*}
By induction on $t$, we have $e(G_{t})>e(T_{n-t,r})+t.$
On the other hand, we have $|V(G_t)|=n-t$, then $e(G_t)\leq \binom{n-t}{2}$. If $n-t<\sqrt{n}$, then $e(G_t)<\binom{\sqrt{n}}{2}=\frac{n-\sqrt{n}}{2}$. Thus
$$n-\sqrt{n}<t\leq e(G_t)\leq\frac{n-\sqrt{n}}{2}$$
 is a contradiction. Therefore, $n-t\geq \sqrt{n}$. For sufficiently large $n$, $|V(G_t)|=n-t\geq \sqrt{n}$ is sufficiently large. Since $|L_t|=0$,
$e(G_t)\leq e(T_{n-t-k+1,r})+(k-1)(n-t)+a-\frac{k}{2}(k-1)$ by Lemma \ref{L is empty in edge extrema}. Recall that $G_{t}$ is obtained from $G_{t-1}$ by deleting a vertex with degree small than $(1-\frac{1}{r}-\varepsilon_2)(n-t+1)$, $G_{t-1}$ is obtained from $G_t$ by adding a new vertex and at most $\big(1-\frac{1}{r}-\varepsilon_2\big)(n-t+1)$ edges incident with the new vertex. Thus
\begin{align*}
e(G_{t-1})&\leq e(G_{t})+\Big(1-\frac{1}{r}-\varepsilon_2\Big)(n-t+1)\\[2mm]
&\leq e(T_{n-t-k+1,r})+(k-1)(n-t)+a-\frac{k}{2}(k-1)+\Big(1-\frac{1}{r}-\varepsilon_2\Big)(n-t+1)\\[2mm]
&< e(T_{n-t-k+2,r})+(k-1)(n-t+1)+a-\frac{k}{2}(k-1).
\end{align*}
By induction on $t$, we have $e(G_0)< e(T_{n-k+1,r})+(k-1)n+a-\frac{k}{2}(k-1)$, and this contradicts the fact $e(G)\geq e(T_{n-k+1,r})+(k-1)n+a-\frac{k}{2}(k-1)$ by Lemma \ref{construct lower bound of edge extrema}. Therefore, $L=\emptyset$. By Lemma \ref{L is empty in edge extrema} we have $G= K_{k-1}\vee G(n-k+1,F)$.
\qed
%
%
%

\section{Proof of Theorem \ref{spectral extrema}}\label{spectral}
In this section, we always suppose that $H$ has the maximum spectral radius among all $n$-vertex $kF$-free graphs. Then we should prove that $e(H)=\mathrm{ex}(n,kF)$  for  sufficiently large $n$.  Recall that in Section \ref{edge}, if $G$ has the maximum number of edges among all $n$-vertex $kF$-free graphs, then $G= K_{k-1}\vee G(n-k+1,F)$, where $G(n-k+1,F)\in \mathrm{EX}(n-k+1,F)$.

Let $\rho(H)$ be the spectral radius of $H$. We first prove that $H$ is connected.

\begin{lemma}\label{connected}
$H$ is connected.
\end{lemma}

\noindent{\bfseries Proof.}
Suppose to the contrary that $H$ is not connected. Assume that $H_1,\ldots,H_s$ are the components of $H$ and $\rho(H_1)=\max\{\rho(H_i):~\ i\in [s]\}$, then $\rho(H)=\rho(H_1)$ and $|V(H_1)|\leq n-1$. For any vertex $u\in V(H_1)$, let $H'$ be the graph obtained from $H_1$ by adding a pendent edge $uv$ at $u$ and $n-1-|V(H_1)|$ isolated vertices. Then $\rho(H')>\rho(H_1)=\rho(H)$ by Lemma \ref{subgraph}. This implies that $H'$ contains a copy of $kF$, say $G_1$, then $G_1$ contains a copy of $F$, say $G_2$, such that $uv\in E(G_2)$. Next, we claim that $d_{H_1}(u)<|V(F)|$. Otherwise, $d_{H_1}(u)\geq |V(F)|$. Then there exists a vertex $w\in N_{H_1}(u)$ and $w\notin V(G_2)$. Let $G_3=G_2\setminus \{uv\}\cup \{uw\}$. Then $G_3\cup (G_1\setminus G_2)$ is a copy of $kF$ in $H_1$, a contradiction. Since $u$ is an arbitrary vertex, $\Delta(H_1)<|V(F)|$.
Thus $\rho(H)=\rho(H_1)\leq \Delta(H_1)<|V(F)|<\rho(K_{k-1}\vee G(n-k+1,F))$. On the other hand, since $H$ has the maximum spectral radius among all $n$-vertex $kF$-free graphs and $K_{k-1}\vee G(n-k+1,F)$ is $kF$-free, then $\rho(H)\geq \rho(K_{k-1}\vee G(n-k+1,F))$, a contradiction. Therefore, $H$ is connected.
\qed

Since $H$ is connected, by the Perron-Frobenius theorem, there exit positive eigenvectors corresponding to $\rho(H)$. Let $\mathbf{x}$ be a positive eigenvector corresponding to $\rho(H)$ with $\max\{x_i: i\in V(H)\}=1$.

\begin{lemma}\label{rho}
Let $H$ be an $n$-vertex $kF$-free graph with maximum spectral radius. Then
$$\rho(H)\geq \frac{r-1}{r}n+O(1).$$
\end{lemma}
\noindent{\bfseries Proof.}
Let $H'$ be an $n$-vertex $kF$-free graph with $e(H')=\mathrm{ex}(n,kF)$.
By Theorem \ref{edge extrema} and (\ref{Rayleigh}),  for sufficiently large $n$,  we have
\begin{align*}
\rho(H)&\geq \rho(H')
\geq \frac{\mathbf{1}^{\mathrm{T}}A(H')\mathbf{1}}{\mathbf{1}^{\mathrm{T}}\mathbf{1}}
= \frac{2e(H')}{n}\\[2mm]
&\geq \frac{2}{n}\left(e(T_{n-k+1,r})+(k-1)n+a-\frac{k}{2}(k-1)\right)\\[2mm]
&\geq\frac{r-1}{r}n+\frac{2(k-1)}{r}-o(n)\\[2mm]
&\geq\frac{r-1}{r}n+O(1).
\end{align*}
\qed

 By Lemmas \ref{stability}, \ref{rho}, Theorem \ref{edge extrema} and a similar argument as in \cite{Kang2023},
 we have the following lemmas.

\begin{lemma}\label{partition}
Let $H$ be an $n$-vertex $kF$-free graph with maximum spectral radius. For every $\gamma >0$
and sufficiently large $n$,
\[
e(H)\geq e(T_{n,r})-\gamma n^2.
\]
Furthermore, there exist a partition $\{V_1, V_2, \cdots, V_{r}\}$ of $V(H)$ and $\gamma_1 >0$   such that $\sum_{1\leq i<j\leq r}e_H(V_i,V_j)$ attains the maximum, and
\[
\sum_{i=1}^{r}e_H(V_i)\leq  \gamma n^2,
\]
and for any $i\in [r]$
$$\frac{n}{r}-\gamma_1 n< |V_i|< \frac{n}{r}+\gamma_1 n.$$
\end{lemma}

\begin{lemma}\label{W}
Suppose that $\gamma$ and $\beta$ are two sufficiently small constants with
$\gamma\leq \beta^2$. Let
\[
Q:=\bigcup_{i=1}^{r}\big\{v\in V_i:~ d_{V_i}(v)\geq 2 \beta n\big\}.
\]
Then $|Q|\leq \beta n$.
\end{lemma}

\begin{lemma}\label{L}
Suppose that $\gamma_2$ is a sufficiently small constant with $\gamma < \gamma_2\ll \beta$. Let
\[
R:=\Big\{v\in V(H):~ d(v)\leq \Big(1-\frac{1}{r}-\gamma_2\Big)n\Big\}.
\]
Then $|R|\leq \gamma_3 n$, where $\gamma_3\ll \gamma_2$ is a sufficiently small constant satisfying $\gamma-\gamma_2 \gamma_3+\frac{r-1}{2r}\gamma_3^2<0$.
\end{lemma}
\noindent{\bfseries Proof.}
Suppose to the contrary that $|R|> \gamma_3 n$, then there exists $R'\subseteq R$ with $|R'|=\lfloor \gamma_3 n \rfloor$. Therefore,
\begin{eqnarray*}
e(H\setminus R')&\geq& e(H)-\sum_{v\in R'}d(v)\\[2mm]
&\geq & e(T_{n,r})-\gamma n^2 - \lfloor \gamma_3 n \rfloor \left(1-\frac{1}{r}-\gamma_2\right)n\\[2mm]
&> & \frac{r-1}{2r}(n-\lfloor \gamma_3 n \rfloor-k+1)^2+(k-1)(n-\lfloor \gamma_3 n \rfloor)+a-\frac{k}{2}(k-1)\\[2mm]
&\geq  & e({K}_{k-1}\vee T_{n-\lfloor \gamma_3 n \rfloor-k+1,r})+a\\[2mm]
&= & \mathrm{ex}(n-\lfloor \gamma_3 n \rfloor,kF).
\end{eqnarray*}
Since $e(H\setminus R')>\mathrm{ex}(n-|R'|,kF)$, $H\setminus R'$ contains a copy of $kF$. This contradicts the fact that $H$ is $kF$-free.
\qed

Let $\ell$ be an integer  that is much larger than $\max\{a,|V(F)|\}$.
Using Lemmas \ref{WL structure}, \ref{WL property} and \ref{WL edge}, we obtain the following  lemmas.

\begin{lemma}\label{spectral WL structure}
$(1)$ For any $u\in Q\setminus R$, $H\setminus (Q\cup R)$ contains a copy of $T_{kr\ell,r}$ such that $u$ is adjacent to all vertices of the copy of $T_{kr\ell,r}$.\\
$(2)$ For $i\in [r]$ and an arbitrary subset $S\subseteq V_{i}\setminus (Q\cup R)$,  $H\setminus (Q\cup R)$ contains a subgraph $H'$ such that $H'$ has a spanning subgraph $T_{r|S|,r}$ with $V(H')\cap V_{i}=S$.
\end{lemma}

\begin{lemma}\label{spectral WL property}
$|Q\setminus R|\leq k-1$, and $d(u)=n-1$ for any $u\in Q\setminus R$.
\end{lemma}

\begin{lemma}\label{independent}
For each $i\in [r]$, $$\Delta(H[V_i\setminus (Q\cup R)])< k\ell,$$ $$\nu(H[V_i\setminus (Q\cup R)])< k\ell,$$ and there exists an independent set $J_i\subseteq V_i\setminus (Q\cup R)$ such that $$|J_i|\geq |V_i\setminus (Q\cup R)|-2(k\ell-1).$$
\end{lemma}

\noindent{\bfseries Proof.}
A similar argument as in the proof of Lemma \ref{WL edge} shows that
 $\Delta(H[V_i\setminus (Q\cup R)])<k\ell$ and  $\nu(H[V_i\setminus (Q\cup R)])< k\ell$ for any $i\in [r]$.
For every $i\in [r]$, let $M^i$ be a maximum matching of $H[V_i\setminus (Q\cup R)]$, and $C^i$ be the set of vertices covered by $M^i$. Then $|C^i|\leq 2(k\ell-1)$.
Let $J_i=V_i\setminus (Q\cup R\cup C^i)$. Then $J_i$ is an independent set and
 $|J_i|\geq |V_i\setminus (Q\cup R)|-2(k\ell-1)$.
\qed

\begin{lemma}\label{HJ}
For any vertex $u\in V(H)$, let $H'$ be the graph with $V(H')=V(H)$ and $E(H')=E(H\setminus \{u\}) \cup \{uw: w\in \cup_{i=2}^{r}J_i\}$. Then $H'$ is $kF$-free.

\end{lemma}
\noindent{\bfseries Proof.}
Otherwise, suppose that $H_1$ is a copy of $kF$ in $H'$. Then $u\in V(H_1)$, and $H_1$ has a copy of $F$, say $H_2$, such that $u\in V(H_2)$. Let $H_3=H_1\setminus H_2$. Then $H_3\cong (k-1)F$ and $u\notin V(H_3)$. Let $N_{H'}(u)\cap V(H_2)=\{w_1,w_2,\ldots, w_s\}$. It is obvious that $w_i\notin V_1$ and $w_i\notin Q\cup R$ for any $i\in [s]$, therefore,
\begin{align*}
d_{V_1}(w_i)
&=d(w_i)-d_{V\setminus V_1}(w_i)\\[2mm]
&\geq \left(1-\frac{1}{r}- \gamma_2 \right)n-2\beta n-(r-2)\left(\frac{n}{r}+\gamma_1 n\right)\\[2mm]
&>\frac{n}{r}-o(n).
\end{align*}
Using Lemma \ref{intersect}, we get
\begin{eqnarray*}
& & \Big|\bigcap_{i=1}^{s}N_{V_1}(w_{i})\setminus (Q\cup R)\Big|\\[2mm]
&\geq &\sum_{i=1}^{s}d_{V_1}(w_{i})-(s-1)|V_1|-|Q|- |R|\\[2mm]
&> & s\left(\frac{n}{r}-o(n)\right)-
(s-1)\left(\frac{n}{r}+\gamma_1 n\right)-\gamma_3 n-(k-1)\\[2mm]
&> & \frac{n}{r}-o(n) > k|V(F)|+1.
\end{eqnarray*}
Thus there exists $u'\in V_1\setminus (Q\cup R)$ and $u'\notin V(H_1)$ such that $u'$ is adjacent to $w_{1},\ldots,w_{s}$. Then $(H_2\setminus \{u\})\cup \{u'\}$ is a copy of $F$ in $H$, and $V((H_2\setminus \{u\})\cup \{u'\})\cap V(H_3)=\emptyset$. Thus $H$ has a copy of $kF$ and this is a contradiction. Therefore $H'$ is $kF$-free.
\qed

\begin{lemma}\label{xv0}
Let $x_{v_0}=\max\{x_v : v\in V(H)\setminus Q\}$. Then $x_{v_0}>1-\frac{2}{r}$ and $v_0\notin R$.

\end{lemma}
\noindent{\bfseries Proof.}
Recall that $\max\{x_v : v\in V(H)\}=1$. Then
\[
\rho(H)\leq |Q|+(n-|Q|)x_{v_0}.
\]
By Lemmas \ref{L} and \ref{spectral WL property}, we have
\begin{eqnarray}\label{Wnumber}
|Q|=|Q\cap R|+|Q\setminus R|\leq |R|+k-1\leq \gamma_3 n+k-1. \label{5.1}
\end{eqnarray}
Combining this with Lemma \ref{rho}, we have
\begin{eqnarray}
x_{v_0}\geq \frac{\rho(H)-|Q|}{n-|Q|}
\geq \frac{\rho(H)-|Q|}{n}
\geq 1-\frac{1}{r}- \gamma_3-\frac{O(1)}{n}
>1-\frac{2}{r}.\label{6.1}
\end{eqnarray}
Therefore, we have
\begin{eqnarray}
\rho (H) x_{v_0}&=&\sum_{v\sim v_0} x_v
= \sum_{\substack{v\in Q, \\ v\sim v_0}} x_v
+ \sum_{\substack{v\notin Q,\\ v\sim v_0} } x_v \nonumber\\[2mm]
&\leq & |Q|+(d(v_0)-|Q|)x_{v_0}. \label{7.1}
\end{eqnarray}

By (\ref{5.1}), (\ref{6.1}), (\ref{7.1}) and Lemma \ref{rho}, we have
\begin{eqnarray*}
d(v_0)&\geq& \rho(H)+|Q|-\frac{|Q|}{x_{v_0}}\\
&\geq& \rho(H)-\frac{2|Q|}{r-2}\\
&\geq& \frac{r-1}{r}n+O(1)-\frac{2\gamma_3 n}{r-2}-\frac{2(k-1)}{r-2}\\
&>& (1-\frac{1}{r}-\gamma_2)n,
\end{eqnarray*}
where the last inequality  holds as $\gamma_3\ll \gamma_2$. Thus $v_0\notin R$.
\qed

\begin{lemma}\label{Jxv0}
$$\sum_{\substack{v\in \cup_{i=2}^{r}J_i } } x_v\geq (\rho (H)-|R|-k\ell-2(k\ell-1)(r-1)) x_{v_0}-|Q|.$$

\end{lemma}
\noindent{\bfseries Proof.}
By Lemma \ref{xv0}, $v_0\in V(H)\setminus (Q\cup R)$. Without loss of generality, we assume that $v_0\in V_1\setminus (Q\cup R)$. By Lemma \ref{independent}, we have
\begin{eqnarray*}
 \rho (H) x_{v_0}
 &=& \sum_{\substack{v\in Q\cup R,\\ v\sim v_0}} x_v  + \sum_{\substack{v\in V_1\setminus (Q\cup R),\\ v\sim v_0}} x_v
 + \sum_{\substack{v\in (\cup_{i=2}^{r}V_i)\setminus (Q\cup R), \\ v\sim v_0} } x_v\\[2mm]
& < & |Q|+|R|x_{v_0}+k\ell x_{v_0}+ \sum_{\substack{v\in (\cup_{i=2}^{r}V_i\setminus J_i)\setminus (Q\cup R), \\ v\sim v_0} } x_v + \sum_{\substack{v\in \cup_{i=2}^{r}J_i, \\ v\sim v_0} } x_v\\[2mm]
& \leq & |Q|+|R|x_{v_0}+k\ell x_{v_0}+ 2(k\ell-1)(r-1)x_{v_0}+ \sum_{\substack{v\in \cup_{i=2}^{r}J_i } } x_v,
\end{eqnarray*}
which implies that
\begin{equation*}
\sum_{\substack{v\in \cup_{i=2}^{r}J_i } } x_v\geq (\rho (H)-|R|-k\ell-2(k\ell-1)(r-1)) x_{v_0}-|Q|.
\end{equation*}
\qed

\begin{lemma}\label{Lemptyset}
$R=\emptyset.$
\end{lemma}
\noindent{\bfseries Proof.}
Suppose to the contrary that there is a vertex $u\in R$, then $d(u)\leq (1-\frac{1}{r}-\gamma_2)n$. Let $H'$ be the graph with $V(H')=V(H)$ and $E(H')=E(H\setminus \{u\}) \cup \{uw: w\in \cup_{i=2}^{r}J_i\}$. Then $H'$ is $kF$-free by Lemma \ref{HJ}.
By Lemmas \ref{rho}, \ref{L}, \ref{xv0} and \ref{Jxv0},
we have
\begin{eqnarray*}
&&\rho(H') - \rho(H)\\[2mm]
&\geq& \frac{\mathbf{x}^T\left(A(H')-A(H)\right)\mathbf{x}}{\mathbf{x}^T\mathbf{x}}\\[2mm]
& =& \frac{2x_{u}}{\mathbf{x}^T\mathbf{x}}\left(\sum_{\substack{w\in \cup_{i=2}^{r}J_i}} x_w - \sum_{v\sim u} x_v\right) \\[2mm]
& \geq& \frac{2x_{u}}{\mathbf{x}^T\mathbf{x}}\Bigl((\rho (H)-|R|-k\ell-2(k\ell-1)(r-1)) x_{v_0}-|Q|-(|Q|+(d(u)-|Q|)x_{v_0}) \Bigr)\\[2mm]
& = &\frac{2x_{u}}{\mathbf{x}^T\mathbf{x}}\Bigl((\rho (H)-|R|-k\ell-2(k\ell-1)(r-1)-d(u)+|Q|) x_{v_0}-2|Q| \Bigr)\\[2mm]
&\geq& \frac{2x_{u}}{\mathbf{x}^T\mathbf{x}}\Bigl( \frac{r-2}{r}(\gamma_2 n- \gamma_3 n- O(1))-2(\gamma_3 n+k-1) \Bigr)>0,
\end{eqnarray*}
where the last inequality holds since $\gamma_3 \ll \gamma_2$. This contradicts the fact that $H$ has the largest spectral radius over all $n$-vertex $kF$-free graphs. So $R$ must be empty.
\qed

\begin{lemma}\label{eigenvector}
For any $v\in V(H)$, $x_v\geq x_{v_0}-\frac{100kr\ell}{n}$.
\end{lemma}
\noindent{\bfseries Proof.}
Suppose to the contrary that there exists $u\in V(H)$ such that $x_u< x_{v_0}-\frac{100kr\ell}{n}$. Let $H'$ be the graph with $V(H')=V(H)$ and $E(H')=E(H\setminus \{u\})\cup \{uw: w\in \cup_{i=2}^{r}J_i\}$. Then $H'$ is $kF$-free by Lemma \ref{HJ}. Since $R=\emptyset$, then $|Q|=|Q\setminus R|\leq k-1$ by Lemma \ref{spectral WL property}. Recall that $x_{v_0}=\max\{x_v : v\in V(H)\setminus Q\}$, by Lemma \ref{Jxv0}, we have
\begin{equation*}
\sum_{\substack{v\in \cup_{i=2}^{r}J_i } } x_v\geq (\rho (H)-k\ell-2(k\ell-1)(r-1)) x_{v_0}-(k-1).
\end{equation*}
 Therefore, we have
\begin{eqnarray*}
&&\rho(H') - \rho(H)\\
& \geq& \frac{2x_{u}}{\mathbf{x}^T\mathbf{x}}\Bigl((\rho (H)-k\ell-2(k\ell-1)(r-1)) x_{v_0}-(k-1)- \rho(H)x_u \Bigr)\\[2mm]
&> &\frac{2x_{u}}{\mathbf{x}^T\mathbf{x}}\Bigl((\rho (H)-k\ell-2(k\ell-1)(r-1)) x_{v_0} -(k-1)- \rho(H)(x_{v_0}-\frac{100kr\ell}{n}) \Bigr)\\[2mm]
&> &\frac{2x_{u}}{\mathbf{x}^T\mathbf{x}}\Bigl( \frac{(r-1)n}{r}\frac{100kr\ell}{n}- k\ell-2(k\ell-1)(r-1)-(k-1)\Bigr)>0.
\end{eqnarray*}
This contradicts the fact that $H$ has the largest spectral radius over all $n$-vertex $kF$-free graphs.
\qed

\begin{lemma}\label{Wk-1}
$|Q|=k-1$.
\end{lemma}
\noindent{\bfseries Proof.}
Let $|Q|=s$. Then $s\leq k-1$ by Lemmas \ref{spectral WL property} and \ref{Lemptyset}. Suppose to the contrary that $s<k-1$.
By Lemmas \ref{independent} and \ref{Lemptyset}, for any $i\in [r]$, $\nu(H[V_i\setminus Q])<k\ell$ and $\Delta(H[V_i\setminus Q])<k\ell$. Thus $\nu(\cup_{i=1}^{r}H[V_i\setminus Q])<rk\ell$ and $\Delta(\cup_{i=1}^{r}H[V_i\setminus Q])<k\ell$. Combining this with Lemma \ref{ffunction}, we have
\begin{align*}
e(\cup_{i=1}^{r}H[V_i\setminus Q])&\leq f(\nu(\cup_{i=1}^{r}H[V_i\setminus Q]),\Delta(\cup_{i=1}^{r}H[V_i\setminus Q]))\\
&< f(kr\ell,k\ell)\\
&< kr\ell(k\ell+1).
\end{align*}
Take $S\subseteq V_1\setminus Q$ with $|S|=k-s-1$. Let $H'$ be the graph with $V(H')=V(H)$ and $E(H')=E(H)\setminus \{uv: uv\in \cup_{i=1}^{r}E(H[V_i\setminus Q])\}\cup \{uv: u\in S,v\in (V_1\setminus Q)\setminus S\}$. It is obvious that $H'$ is $kF$-free. Therefore,
\begin{align*}
\rho(H') - \rho(H)
&\geq \frac{\mathbf{x}^T\left(A(H')-A(H)\right)\mathbf{x}}{\mathbf{x}^T\mathbf{x}}\\[2mm]
&=\frac{2}{\mathbf{x}^T\mathbf{x}}\left(\sum_{ij\in E(H')}x_ix_j-\sum_{ij\in E(H)}x_ix_j\right) \\[2mm]
& \geq \frac{2}{\mathbf{x}^T\mathbf{x}}\left(|S|(|V_1|-|Q|-|S|)(x_{v_0}-\frac{100kr\ell}{n})^2- kr\ell(k\ell+1)\right) \\[2mm]
& \geq \frac{2}{\mathbf{x}^T\mathbf{x}}\Bigl( (k-s-1)(\frac{n}{r}-o(n))(\frac{r-2}{r}-o(1))^2-kr\ell(k\ell+1)\Bigr)\\[2mm]
&>0.
\end{align*}
This contradicts the fact that $H$ has the largest spectral radius over all $n$-vertex $kF$-free graphs.
Therefore, $|Q|=s=k-1$.
\qed

\begin{lemma}\label{Bi}
For any $i\in [r]$, let $T_i=\{u\in V_i\setminus Q: d_{V_i\setminus Q}(u)=0\}$. Then $u$ is adjacent to all vertices of $V(H)\setminus V_i$ for every vertex $u\in T_i$.
\end{lemma}
\noindent{\bfseries Proof.} Since each vertex of $Q$ has degree $n-1$, we only need to prove that each vertex of $T_i$ is adjacent to all vertices of $(V(H)\setminus Q)\setminus V_i$ for any $i\in [r]$. Suppose to the contrary that there exist $i_0\in [r]$ and $v\in T_{i_0}$ such that there is a vertex $w_{1}\notin V_{i_0}\cup Q$ and $vw_{1}\notin E(H)$. Without loss of generality, we may assume that $i_0=1$. Let $H'$ be the graph with $V(H')=V(H)$ and $E(H')=E(H)\cup \{vw_{1}\}$.
Now we prove that $H'$ is $kF$-free. Otherwise, $H'$ contains a  copy of $kF$, say $H_1$. Then $H_1$ has a copy of  $F$, say $H_2$, such that $vw_{1}\in E(H_2)$. Let $H_3=H_1\setminus H_2$. Then $H_3$ is a copy of $(k-1)F$ in $H$. We first claim that $V(H_2)\cap Q=\emptyset$. Otherwise, by Lemma \ref{spectral WL structure}, we can find a copy of $F$ in $H$, say $H_4$, such that $V(H_4)\cap V(H_3)=\emptyset$ and $H_3\cup H_4$ is a copy of $kF$ in $H$, a contradiction.
Let $N_{H'}(v)\cap V(H_2)=\{w_1,w_2,\ldots, w_s\}$. It is obvious that $w_i\notin V_1$ and $w_i\notin Q$ for any $i\in [s]$.
Using the similar discussion as in Lemma \ref{HJ}, there exists $v'\in V_1\setminus Q$ and $v'\notin V(H_1)$ such that $v'$ is adjacent to $w_{1},\ldots,w_{s}$. Then $(H_2\setminus \{v\})\cup \{v'\}$ contains a copy of $F$, and $V((H_1\setminus \{v\})\cup \{v'\})\cap V(H_3)=\emptyset$. Thus $H$ has a copy of $kF$, this is a contradiction. Thus $H'$ is $kF$-free. From the construction of $H'$, we see that $\rho(H')>\rho(H)$, which contradicts the assumption that $H$ has the maximum spectral radius among all $n$-vertex $kF$-free graphs.
\qed

By Lemmas \ref{spectral WL property}, \ref{Lemptyset} and \ref{Wk-1}, $H[Q]= K_{k-1}$. Let $H_{in}=\cup_{i=1}^{r}H[V_i\setminus Q]$, $K$ be the graph with $V(K)=V(H)\setminus Q$ and $E(K)=\cup_{1\leq i<j\leq r}\{uv:~ u\in V_i\setminus Q, v\in V_j\setminus Q\}$.

\begin{lemma}\label{Hout}
Let $H_{out}$ be the graph with $V(H_{out})=V(K)$ and $E(H_{out})=E(K)\setminus E(H\setminus Q)$. Then $e(H_{out})\leq 2a^2r^2$ and $e(H_{in})-e(H_{out})\leq a$.

\end{lemma}
\noindent{\bfseries Proof.}
We first prove the following claim.

\noindent{\bfseries Claim A.}
For any $i\in [r]$, $e(H[V_i\setminus Q])\leq a$.

\noindent{\bfseries Proof of Claim A. }
Suppose to the contrary that there is an $i_0\in [r]$ such that $e(H[V_{i_0}\setminus Q])>a$. Without loss of generality, we may assume that $e(H[V_1\setminus Q])> a.$
By Lemmas \ref{partition} and \ref{Wk-1}, we have $|V_i\setminus Q|\geq \left(\frac{1}{r}-  \varepsilon_1\right)n -k+1 \ge  \ell$ for any $i \in [r]$. Let $u_{1,1},u_{1,2}, \ldots, u_{1,\ell} $ be $\ell$ vertices chosen  from $V_1 \setminus Q$ such that the induced subgraph of $\{u_{1,1},u_{1,2}, \ldots, u_{1,\ell}\} $ in $H$ contains at least $a+1$ edges. Using Lemma \ref{spectral WL structure} (2) for $S=\{u_{1,1},u_{1,2}, \ldots, u_{1,\ell}\}$,
we obtain a subgraph $H'$ such that $H'$ has a spanning subgraph $T_{r\ell,r}$, and $e(H')> e(T_{r\ell,r})+a$. This implies that $H'$ contains a copy of $F$, say $H_1$. Since $|Q|=k-1$ and each vertex of $Q$ has degree $n-1$, we can find a copy of $(k-1)F$, say $H_2$, such that $V(H_1)\cap V(H_2)=\emptyset$. Thus $H_1\cup H_2$ is a copy of $kF$ in $H$, a contradiction. Therefore, for each $i \in [r]$, $e(H[V_i\setminus Q])\leq a$.
\qed

For any $i\in [r]$, let $B_i=(V_i\setminus Q)\setminus T_i$. Then by Claim A, we have
\[
|B_i|\leq  \sum_{u\in B_i}d_{V_i}(u)=\sum_{u\in V_i\setminus Q}d_{V_i}(u)=2e(H[V_i\setminus Q])<2a.
\]
By Lemma \ref{Bi}, each vertex of $T_i$ is adjacent to all vertices of $V(H)\setminus V_i$ for any $i\in [r]$. Thus $$e(H_{out})\leq \sum_{1\leq i<j\leq r}|B_i||B_j|< \binom{r}{2}(2a)^2\leq 2a^2r^2.$$

Now we prove that $e(H_{in})-e(H_{out})\leq a$. Suppose to the contrary that $e(H_{in})-e(H_{out})> a.$ For each $i\in [r]$, let $S_i$ be the vertex set satisfying
$B_i\subseteq S_i\subseteq V_i\setminus Q$ and $|S_i|=s$. Let $S=\cup_{i=1}^{r}S_i$, $H'=H[S]$. By Lemma \ref{Bi}, we have
$e(H')\geq  e(T_{rs,r})+e(H_{in})-e(H_{out})>e(T_{rs,r})+a$. Thus $H'$ contains a copy of $F$, say $H_2$. Since $|Q|=k-1$ and each vertex of $Q$ has degree $n-1$, we can find a copy of $(k-1)F$ such that this $(k-1)F$ is disjoint with $H_2$, that is $H$ has a copy of $kF$, a contradiction.
\qed

\begin{lemma}\label{balance}
For any $1\leq i<j\leq r$,  $\left||V_i\setminus Q|-|V_j\setminus Q|\right|\leq 1$.
\end{lemma}

\noindent{\bfseries Proof.}
For any $i\in [r]$, let $|V_i\setminus Q|=n_i$. Then by the definition of $K$, $K$ is a complete $r$-partite graph with classes of sizes $n_1,\cdots, n_r$, denoted by $K_{r}(n_1,\cdots, n_r)$.
Suppose that $n_1\geq n_2\geq \cdots \geq n_{r}$.  Now we prove the result by contradiction.
Assume that there exist $i_0, j_0$ with $1\leq i_0 < j_0\leq r$ such that $n_{i_0}-n_{j_0}\geq 2$. Let $K'=K_{r}(n_1,\cdots, n_{i_0}-1,\cdots,n_{j_0}+1,\cdots,n_{r})= K_{r}(n'_1,n'_2,\cdots,n'_{r})$, where $n'_1\geq n'_2\geq \cdots \geq n'_{r}$.

\noindent{\bfseries Claim 1.} There exists a constant $c_1>0$ such that $$\rho(K_{k-1}\vee T_{n-k+1,r})-\rho(K_{k-1}\vee K)\geq \frac{c_1}{n}.$$

\noindent{\bfseries Proof of Claim 1.}
Let $\mathbf{y}$ be a positive  eigenvector of $K_{k-1}\vee K$ corresponding to $\rho(K_{k-1}\vee K)$. By the symmetry we may assume $\mathbf{y}=(\underbrace{y_1,\cdots,y_1}_{n_1},\underbrace{y_2,\cdots,y_2}_{n_2},\cdots, \underbrace{y_r,\cdots,y_r}_{n_r},\underbrace{y_{r+1},\cdots,y_{r+1}}_{k-1})^{\mathrm{T}}$.
By  (\ref{eigenequation}), we have
\begin{eqnarray*}
\rho(K_{k-1}\vee K) y_i=\sum_{j=1}^{r}n_jy_j-n_iy_i+(k-1)y_{r+1}, \text{ for any } i\in [r],\\ 
\rho(K_{k-1}\vee K) y_{r+1}=\sum_{j=1}^{r}n_jy_j+(k-2)y_{r+1}.
\end{eqnarray*}
Thus $y_i=\frac{\rho(K_{k-1}\vee K)+1}{\rho(K_{k-1}\vee K)+n_i}y_{r+1}$ for any $i\in [r]$, which implies that $y_{r+1}=\max\{y_v :~ v\in V(K_{k-1}\vee K)\}$. Without loss of generality, we assume that $y_{r+1}=1$, then
\begin{eqnarray}
y_i=\frac{\rho(K_{k-1}\vee K)+1}{\rho(K_{k-1}\vee K)+n_i}, \text{ for any } i\in [r]. \label{116}
\end{eqnarray}
Let $u_{i_0}\in V_{i_0}\setminus Q$ be a fixed vertex. Then $K_{k-1}\vee K'$ can be obtained from $K_{k-1}\vee K$ by deleting all edges between $u_{i_0}$ and $V_{j_0}\setminus Q$, and adding all edges between $u_{i_0}$ and $V_{i_0}\setminus (Q\cup\{u_{i_0}\})$.
According to  (\ref{Rayleigh}), we deduce that
\begin{eqnarray*}
&&\rho(K_{k-1}\vee K')-\rho(K_{k-1}\vee K)\\[2mm]
&\geq& \frac{2}{\mathbf{y}^{\mathrm{T}}\mathbf{y}}\left((n_{i_0}-1)y_{i_0}^2-n_{j_0}y_{i_0}y_{j_0}\right)\\[2mm]
&=& \frac{2y_{i_0}}{\mathbf{y}^{\mathrm{T}}\mathbf{y}}\left((n_{i_0}-1)\frac{\rho(K_{k-1}\vee K)+1}{\rho(K_{k-1}\vee K)+n_{i_0}}-n_{j_0}\frac{\rho(K_{k-1}\vee K)+1}{\rho(K_{k-1}\vee K)+n_{j_0}}\right)\\[2mm]
&\geq& \frac{2y_{i_0}}{\mathbf{y}^{\mathrm{T}}\mathbf{y}}\frac{(\rho(K_{k-1}\vee K)+1)(\rho(K_{k-1}\vee K)-n_{j_0})}{(\rho(K_{k-1}\vee K)+n_{i_0})(\rho(K_{k-1}\vee K)+n_{j_0})},
\end{eqnarray*}
where the last inequality holds as $n_{i_0}-n_{j_0}\geq 2$. Since $\delta(K_{k-1}\vee K)\leq \rho(K_{k-1}\vee K)\leq \Delta(K_{k-1}\vee K)$,  $\rho(K_{k-1}\vee K)=\Theta(n)$.
  From (\ref{116}) and the fact that $\mathbf{y}^{\mathrm{T}}\mathbf{y}\leq n$,  it follows that there exists   a constant $c_1>0$ such that $\rho(K_{k-1}\vee K')-\rho(K_{k-1}\vee K)\geq \frac{c_1}{n}$. Therefore, $\rho(K_{k-1}\vee T_{n-k+1,r})-\rho(K_{k-1}\vee K)\geq\rho(K_{k-1}\vee K')-\rho(K_{k-1}\vee K)\geq \frac{c_1}{n}$.\qed

\noindent{\bfseries Claim 2.} There exists a constant $c_2>0$ such that $$\rho(K_{k-1}\vee T_{n-k+1,r})-\rho(K_{k-1}\vee K))\leq \frac{c_2}{n^2}.$$

\noindent{\bfseries Proof of Claim 2.}
According to the definitions of $H_{in}$, $H_{out}$ and $K$, we have $E(H)=E(K_{k-1}\vee K)\cup E(H_{in})\setminus E(H_{out})$. Recall that $x_{v_0}=\max\{x_v : v\in V(H)\setminus Q\}$.
Combining with Lemmas \ref{eigenvector} and \ref{Hout}, we have 
\begin{align}
\rho(H)-\rho(K_{k-1}\vee K) 
&\leq \frac{2\sum_{ij\in E(H_{in})}\mathbf{x}_i\mathbf{x}_j}{\mathbf{x}^{\mathrm{T}}\mathbf{x}}- \frac{2\sum_{ij\in E(H_{out})}\mathbf{x}_i\mathbf{x}_j}{\mathbf{x}^{\mathrm{T}}\mathbf{x}} \nonumber\\[2mm]
&\leq  \frac{2e(H_{in})x_{v_0}^2}{\mathbf{x}^{\mathrm{T}}\mathbf{x}}- \frac{2e(H_{out})(x_{v_0}-\frac{100kr\ell}{n})^2}{\mathbf{x}^{\mathrm{T}}\mathbf{x}} \nonumber\\[2mm]
&\leq \frac{2(e(H_{in})-e(H_{out}))x_{v_0}^2}{n(x_{v_0}-\frac{100kr\ell}{n})^2}+\frac{2e(H_{out})\frac{200kr\ell}{n}}{n(x_{v_0}-\frac{100kr\ell}{n})^2}\nonumber\\[2mm]
&\leq \frac{2a}{n-\frac{200kr\ell}{x_{v_0}}}+\frac{800r^3a^2k\ell}{n^2x_{v_0}^2-200kr\ell n x_{v_0}}\label{GC}
\end{align}
On the other hand, let $\mathbf{z}$ be an eigenvector corresponding to $\rho(K_{k-1}\vee T_{n-k+1,r})$ with $\max\{z_v : v\in V(K_{r-1}\vee T_{n-k+1,r})\}=1$. Let  $m_1=\lfloor\frac{n-k+1}{r}\rfloor$, $m_2=\lceil\frac{n-k+1}{r}\rceil$, and $b=n-k+1-rm_1$. By the symmetry we may assume $\mathbf{z}=(\underbrace{z_1,\cdots,z_1}_{(r-b)m_1},\underbrace{z_2,\cdots,z_2}_{bm_2}, \underbrace{z_{3},\cdots,z_{3}}_{k-1})^{\mathrm{T}}$. By the similar
discussion as in Claim 1,
 we have $z_3=1$,
\[
z_1= \frac{\rho(K_{k-1}\vee T_{n-k+1,r})+1}{\rho(K_{k-1}\vee T_{n-k+1,r})+m_1}>
z_2= \frac{\rho(K_{k-1}\vee T_{n-k+1,r})+1}{\rho(K_{k-1}\vee T_{n-k+1,r})+m_2}> 1-\frac{1}{r}.
\]
Let $H'\in \mathrm{EX}(n,kF)$. Then $H'\cong K_{k-1}\vee G(n-k+1,F)$, where $G(n-k+1,F)\in \mathrm{EX}(n-k+1,F)$.
Using the similar argument as the proof of Theorem 1.3 in  \cite{Kang2023},
  we have 
\begin{align*}
\rho(H)\geq \rho(H')
\geq \frac{\mathbf{z}^\mathrm{T}(A(K_{k-1}\vee T_{n-k+1,r}))\mathbf{z}}{\mathbf{z}^{\mathrm{T}}\mathbf{z}}+\frac{2az_2^2}{\mathbf{z}^{\mathrm{T}}\mathbf{z}}= \rho(K_{k-1}\vee T_{n-k+1,r})+\frac{2az_2^2}{\mathbf{z}^{\mathrm{T}}\mathbf{z}}.
\end{align*}
Since $m_2-m_1\leq 1$, $bm_2+(r-b)m_1=n-k+1$ and $\rho(K_{k-1}\vee T_{n-k+1,r})\geq \delta(K_{k-1}\vee T_{n-k+1,r})=n-m_2$, we have
\begin{align*}
\frac{\mathbf{z}^{\mathrm{T}}\mathbf{z}}{z_2^2}&\leq \frac{k-1}{z_2^2}+ \frac{(r-b)m_1z_1^2}{z_2^2}+ bm_2\\[2mm]
&=n+(k-1)\left(\frac{1}{z_2^2}-1\right)+(r-b)m_1\left(\frac{z_1^2}{z_2^2}-1\right)\\
&\leq n+(k-1)\left((\frac{r}{r-1})^2-1\right)+(r-b)m_1\left((1+\frac{1}{n-1})^2-1\right)\\
&\leq n+ c_3,
\end{align*}
where $c_3>0$ is a constant. Therefore, for sufficiently large $n$, there exists a constant $c_2>0$ such that
\begin{equation}\label{Tnr-1}
\rho(H)-\rho(K_{k-1}\vee T_{n-k+1,r})\geq \frac{2a}{n+ c_3}\geq \frac{2a}{n}\left(1-\frac{c_2}{n}\right).
\end{equation}

It follows from (\ref{GC}), (\ref{Tnr-1}) and $x_{v_0}>1-\frac{2}{r}$ that there exists a constant $c_4>0$ such that
\begin{eqnarray*}
& &\rho(K_{k-1}\vee T_{n-k+1,r})-\rho(K_{k-1}\vee K)\\[2mm]
&\leq &\frac{2a}{n-\frac{200kr\ell}{x_{v_0}}}+\frac{800r^3a^2k\ell}{n^2x_{v_0}^2-200kr\ell n x_{v_0}}-\frac{2a}{n}\left(1-\frac{c_2}{n}\right)\\[2mm]
&< & \frac{2a\frac{200kr^2\ell}{r-2}}{n(n-\frac{200kr^2\ell}{r-2})}+\frac{800r^3a^2k\ell}{\frac{(r-2)^2n^2}{r^2}-200kr\ell n}+\frac{2ac_2}{n^2}\\[2mm]
&\leq & \frac{c_4}{n^2}.
\end{eqnarray*}
\qed

Combining Claim 1 and Claim 2, we have
 \begin{align*}
 \frac{c_1}{n}\leq \rho(K_{k-1}\vee T_{n-k+1,r})-\rho(K_{k-1}\vee K) \leq \frac{c_4}{n^2},
 \end{align*}
 which is a contradiction when $n$ is sufficiently large. Thus $\left||V_i\setminus W|-|V_j\setminus W|\right|\leq 1$ for any $1\leq i<j\leq r$.
\qed

\medskip
\noindent{\bfseries Proof of Theorem \ref{spectral extrema}.} Now we  prove that $e(H)=\mathrm{ex}(n,kF)$. Otherwise, we assume that $e(H)\leq \mathrm{ex}(n,kF)-1$. Let $H'$ be a $kF$-free graph with $e(H')=\mathrm{ex}(n,kF)$. Then $H'= K_{k-1}\vee G(n-k+1,F)$, where $G(n-k+1,F)\in \mathrm{EX}(n-k+1,F)$. Using the similar argument as the proof of Theorem 1.3 in  \cite{Kang2023},
 we get that
 $H'$ can be obtained from $K_{k-1}\vee T_{n-k+1,r}$ by adding at most $ar$ edges and deleting at most $a(r-1)$ edges
 in $T_{n-k+1,r}$. By Lemma \ref{balance}, we may assume that $\{V_1, V_2,\cdots, V_{r}\}$ is a  partition of $V(H')$. Let $E_1=E(H)\setminus E(H')$, $E_2=E(H')\setminus E(H)$. Then $E(H')=(E(H)\cup E_2)\setminus E_1$, and
\[
|E(H)\cap E(H')|+|E_1|=e(H)<e(H')=|E(H)\cap E(H')|+|E_2|,
\]
which implies that
$|E_2|\geq |E_1|+1$. Furthermore, by Lemma \ref{Hout}, we have
\begin{eqnarray*}
|E_2|=e(H')-e(H)\leq ar+e(H_{out})< ar+2a^2r^2< 3a^2r^2.
\end{eqnarray*}
Therefore,
\begin{align*}
\rho(H')-\rho(H)
&\geq \frac{2\sum_{ij\in E_2}\mathbf{x}_i\mathbf{x}_j}{\mathbf{x}^{\mathrm{T}}\mathbf{x}}-\frac{2\sum_{ij\in E_1}\mathbf{x}_i\mathbf{x}_j}{\mathbf{x}^{\mathrm{T}}\mathbf{x}}\\[2mm]
&\geq \frac{2}{\mathbf{x}^{\mathrm{T}}\mathbf{x}}\Big(|E_2|(x_{v_0}-\frac{100kr\ell}{n})^2- |E_1|x^2_{v_0}\Big)\\[2mm]
&\geq  \frac{2}{\mathbf{x}^{\mathrm{T}}\mathbf{x}}\Big(|E_2|x^2_{v_0}-\frac{200kr\ell x_{v_0}}{n}|E_2|- |E_1|x^2_{v_0}\Big)\\[2mm]
&\geq \frac{2x_{v_0}}{\mathbf{x}^{\mathrm{T}}\mathbf{x}}\Big(x^2_{v_0}-\frac{200kr\ell}{n}|E_2|\Big)\\[2mm]
&\geq \frac{2x_{v_0}}{\mathbf{x}^{\mathrm{T}}\mathbf{x}}\Big((1-\frac{2}{r})^2-\frac{200kr\ell}{n}  3a^2r^2\Big)\\[2mm]
&> 0,
\end{align*}
which contradicts the assumption that $H$ has the maximum spectral radius among all $n$-vertex $kF$-free graphs. Hence $e(H)=\mathrm{ex}(n,kF)$.
\qed

\end{document}